\documentclass{amsart}
\usepackage{amssymb}
\begin{document}
\newtheorem{theo}{Theorem}[section]
\newtheorem{prop}[theo]{Proposition}
\newtheorem{lemma}[theo]{Lemma}
\newtheorem{exam}[theo]{Example}
\newtheorem{coro}[theo]{Corollary}
\theoremstyle{definition}
\newtheorem{defi}[theo]{Definition}
\newtheorem{rem}[theo]{Remark}

\newcommand{\en}{{\mathbb{N}}}
\newcommand{\ze}{{\mathbb{Z}}}
\newcommand{\A}{{\mathcal A}}
\newcommand{\C}{{\mathcal C}}
\newcommand{\es}{{\mathcal S}}
\newcommand{\re}{{\mathcal R}}
\newcommand{\et}{{\mathcal T}}
\newcommand{\Be}{{\mathcal B}}
\newcommand{\Ha}{{\mathcal H}}
\newcommand{\Ka}{{\mathcal K}}
\newcommand{\pf}{{\it Proof:\/ }}
\newcommand{\cspa}{\mathop{\overline {\rm span}\ }}
\newcommand{\Mcc}{{M\raise .45ex\hbox{c}}}
\newcommand{\zuz}{\upharpoonright}

\author{David Op\v ela}

\title [And\^o's Theorem and Parrott's Example]
       {A Generalization of And\^o's Theorem \\ and Parrott's Example}

\address {Department of Mathematics\\
          Campus Box 1146\\
          Washington University in Saint Louis\\
          Saint Louis, MO 63130\\
          U.S.A.}

\email {opela@math.wustl.edu}

\date{\em December 17, 2004}



\begin{abstract}
And\^o's theorem states that any pair of commuting contractions on a Hilbert
space can be dilated to a pair of commuting unitaries. Parrott presented an
example showing that an analogous result does not hold for a triple of
pairwise commuting contractions. We generalize both of these results as
follows. Any $n$-tuple of contractions that commute according to a graph
without a cycle can be dilated to an $n$-tuple of unitaries that commute
according to that graph. Conversely, if the graph contains a cycle, we
construct a counterexample.
\end{abstract}

\maketitle

\section {Introduction}

Foia\c s and Sz.-Nagy's theory models Hilbert space operators as ``parts'' of
simpler operators. We will begin by recalling its main results, for the proofs
see Chapter 10 in \cite {AM}. Very detailed treatment is in \cite {SF}. Recent
books \cite {Pau}, \cite {Pi2} contain up-to-date expositions of some aspects
of the theory.

If $\Ha \subset \Ka$ are two Hilbert spaces, $T \in \Be(\Ha)$, $W \in \Be(\Ka)$
operators, we say that $W$ is an {\it extension\/} of $T$, if $T$ is the
restriction of $W$ to $\Ha$, i.e., if $Tx = Wx$ for all $x \in \Ha$. With the
same notation, we say that $W$ is a {\it dilation\/} of $T$, or that $T$ is
a {\it compression} of $W$, if $T^n = PW^n \zuz \Ha$, for all $n \ge 0$, where
$P$ is the orthogonal projection from $\Ka$ onto $\Ha$. By a result of Sarason,
we can say equivalently, that $W$ has the following structure 
$$ W = \begin{pmatrix}     T      &  0    &  \ast  \\
                           \ast   &  \ast &  \ast  \\
                           0      &  0    &  \ast
       \end{pmatrix}
       \begin{matrix}      \dots  &  \Ha  \\
                           \dots  &  \tilde \Ka \ominus \Ha \\
                           \dots  &  \Ka \ominus \tilde \Ka
       \end{matrix},$$
that is, $\Ha$ is the orthogonal difference of two invariant subspaces for $W$.

Recall that $A \in \Be(\Ha)$ is a {\it co-isometry\/} if $A^\ast$ is an
isometry, that is, if $AA^\ast = I$.

The first result of the theory, due to Sz.-Nagy, asserts that every contraction
$T \in \Be(\Ha)$ has a co-isometric extension. In fact, there is a
{\it minimal\/} co-isometric extension $W_0 \in \Be(\Ka_0)$ of $T$ that is
characterized by 
$$\Ka_0 = \cspa \{(W_0^\ast)^n \Ha;\ n \in \en\}.$$
Any two minimal co-isometric extensions $W_0$ and $\widetilde W_0$ are
unitarily equivalent via a unitary that restricts to the identity on $\Ha$ ---
we denote this by $W_0 \cong \widetilde W_0$. For any co-isometric extension
$W$ of $T$ we have $W \cong W_0 \oplus \widetilde W$, where $\widetilde W$ is
another co-isometry. If $T$ is already a co-isometry, then any of its
co-isometric extensions is the direct sum of $T$ and another co-isometry. A
minimal co-isometric extension of an isometry is a unitary. If $T$ is a
contraction, denote by $V$ any of the co-isometric extensions of the
contraction $T^\ast$. Let $U$ be a minimal co-isometric extension of the
isometry $V^\ast$, then one easily checks (by Sarason's characterization) that
$U$ is a dilation of $A$. The Sz.-Nagy dilation theorem follows --- any
contraction has a unitary dilation. Given $T \in \Be(\Ha)$ there is a
{\it minimal\/} unitary dilation $U_0$, unique up to a unitary that restricts
to $I_\Ha$, such that any unitary dilation $U$ of $T$ has the form $U \cong U_0
\oplus U_1$.
  
Let us stress that both the extension and the dilation result above say that
a general contraction is a part of a simpler operator. Indeed, unitaries are
well-understood via the spectral theorem. As for co-isometries, by the Wold
decomposition, any co-isometry is the orthogonal sum of (a finite or infinite
number of) copies of the unilateral backward shift and a unitary.

Both of the above results can be generalized to a pair of commuting
contractions. And\^o's theorem \cite {A} states that given a pair of commuting
contractions, we can extend each of them to a co-isometry such that the
two co-isometries commute. As in the single operator case, if the contractions
are isometries, the co-isometries can be constructed to be unitaries. The
dilation version asserts that for any pair of commuting contractions $A, B$, we
can find commuting unitaries $U, V$ such that
$$ PU^mV^n \zuz \Ha = A^mB^n, \hbox { for all } m, n \ge 0.$$
Note that the equality above implies that $U$ (resp. $V$) is a dilation of $A$
(resp. $B$) by taking $n=0$ ($m=0$, respectively). However, not every pair of
dilations satisfies this property. This more restrictive relation is more
desirable, since it implies that the map $U \mapsto A, V \mapsto B$ extends to
an algebra homomorphism between the operator algebras generated by $U, V$ and
$A, B$, respectively.

There is a related theorem whose modification we will use. The commutant
lifting theorem of Foia\c s and Sz.-Nagy asserts that given a pair of commuting
contractions and a co-isometric extension (or unitary dilation) of one of them,
one can extend (or dilate) the other one to a contraction that commutes with the
given co-isometric extension (unitary dilation, respectively).

Surprisingly, And\^o's theorem cannot be generalized to three (or more)
commuting contractions. The first counterexample was constructed by S.~Parrott,
see \cite {Par}. There are some sufficient conditions on when an $n$-tuple of
commuting contractions can be dilated to an $n$-tuple of commuting unitaries,
e.g., if the operators doubly commute, see Theorem 12.10 in \cite {Pau}.

There are generalizations of And\^o's result to an $n$-tuple of contractions.
Ga\c spar and R\'acz assume only that the $n$-tuple is {\it cyclic
commutative}, see \cite {GR}. Their result was further generalized by
G.~Popescu \cite {Po}. Instead of starting with $n$-tuples of contractions, one
can work with the so-called {\it row contractions\/}, that is, with $n$-tuples
satisfying $\sum_j T_jT_j^\ast \le I$. This case has been extensively studied
--- see the recent survey \cite {B} and the references therein.

We derive a different generalization of And\^o's theorem, namely, we assume
that only some of the $\binom n2$ pairs commute --- see the definition below.

Before we will be able to formulate our results, let us recall a few basic
definitions and facts from graph theory. A {\it graph\/} $G$ is a pair $(V(G),
E(G))$, where $V(G)$ is a set (the elements of $V(G)$ are called
{\it vertices\/} of $G$) and $E(G)$ is a set of unordered pairs of distinct
vertices --- these are called {\it edges\/}. A graph $G'$ is a {\it subgraph\/}
of a graph $G$, if $V(G') \subset V(G)$ and $E(G') \subset E(G)$. A {\it cycle
of length $n$\/} is the graph $G$ with $V(G) = \{v_1, \dots, v_n\}$ and $E(G) =
\{(v_j, v_{j+1});\ 1 \le j \le n-1\} \cup \{(v_1, v_n)\}$. A graph $G$ is 
{\it connected,\/} if for any two vertices $v, w$ there exists a sequence of
vertices $\{v_k\}_{k=0}^m \subset V(G)$ with $v_0 = v$, $v_m = w$ and such that
$(v_j, v_{j+1}) \in E(G)$ for $j = 0,\dots, m-1$. A graph is {\it acyclic\/},
if it does not contain a cycle as a subgraph, and a connected acyclic graph is
called a {\it tree\/}. Every tree has a vertex that lies on exactly one edge.

\begin{defi}
Let $A_1, A_2, \dots, A_n \in \Be(\Ha)$ be an $n$-tuple of operators, and let
$G$ be a graph on the vertices $\{1, 2, \dots, n\}$. We say that the operators
$A_1, \dots, A_n$ {\it commute according to $G$,}\/ if $A_iA_j = A_jA_i$
whenever $(i,j)$ is an edge of $G$.
\end{defi}

We prove that given $G$, every $n$-tuple of contractions commuting according to
$G$ has unitary dilation that commute according to $G$, if and only if $G$ is
acyclic. And\^o's theorem is a special case, when $G$ is the acyclic graph
consisting of two vertices and an edge connecting them. The case of three
commuting contractions, Parrott's example, corresponds to the cycle of length
three. If $G$ is a graph with no edges, the result is also known, see Exercise
5.4, p.~71 in \cite {Pau}.

\section {The Main Result}

Our main result is Theorem \ref {main}. We will make use of the following lemma
that can be regarded as a version of the commutant lifting theorem.

\begin{lemma} \label {wcl}
Let $A, B \in \Be(\Ha)$ be commuting contractions and let $\tilde X \in
\Be(\tilde \Ka)$ be a co-isometric extension of $A$. Then there exists a
Hilbert space $\Ka$ containing $\tilde \Ka$ and commuting co-isometries $X, Y
\in \Be(\Ka)$ such that $X$ extends $\tilde X$ and $Y$ extends $B$. Moreover,
if $A, B$ are isometries and $\tilde X$ is a unitary, we can construct $X$ and
$Y$ to be unitary.
\end{lemma}

\pf Let $X_0 \in \Be(\Ka_0)$ be a minimal co-isometric extension of $A$, then
$\tilde X \cong X_0 \oplus X_1$ where $X_1 \in \Be(\Ka_1)$ is a co-isometry.
By taking a different minimal co-isometric extension, we may assume that
$\tilde X = X_0 \oplus X_1$. By And\^o's theorem there exist commuting
co-isometries $X_A, Y_A \in \Be(\Ka_A)$ extending $A, B$, respectively. Hence,
$X_A = X_0 \oplus X_2$, where $X_2 \in \Be(\Ka_2)$ is a co-isometry and
$Y_A$ has a corresponding decomposition
$$Y_A = \begin {pmatrix} Y_{11} & Y_{12} \\
                         Y_{21} & Y_{22}
        \end {pmatrix}.
$$
We set $\Ka := \Ka_0 \oplus \Ka_2 \oplus \Ka_1$,
$$ X := X_0 \oplus X_2 \oplus X_1,
   \hskip 1cm \hbox {and} \hskip 1cm
   Y :=
   \begin {pmatrix} Y_{11} & Y_{12} &         0 \\
                    Y_{21} & Y_{22} &         0 \\
                         0 &      0 & I_{\Ka_1}
   \end {pmatrix}.
$$
Clearly, $X, Y$ extend $\tilde X, B$, respectively. Also, it is easy to check
that $X$ and $Y$ are commuting co-isometries.

To prove the second part, one follows the above proof, uses the fact that $X_A,
Y_A$ can be chosen to be unitary and then verifies that $X, Y$ will also be
unitary. \qed
\medskip

We will construct the unitary dilations using co-isometric extensions in a way
sketched in the introduction for a single contraction. Thus we need the
following lemma.

\begin{lemma} \label {ind}
Let $G$ be a graph without a cycle on vertices $\{1, 2, \dots, n\}$ and let
$A_1, \dots, A_n \in \Be(\Ha)$ an $n$-tuple of contractions that commute
according to $G$. Then there exist a Hilbert space $\Ka$ containing $\Ha$ and
an $n$-tuple of co-isometries $X_1, \dots, X_n \in \Be(\Ka)$ that commute
according to $G$ and such that $X_j$ extends $A_j$, for $j = 1, 2, \dots, n$.
Moreover, if the $A_j$'s are all isometries, the $X_j$'s can be chosen to be
unitaries.
\end{lemma}

\pf We may assume that $G$ is connected since the result for trees easily
implies the general case. Indeed, otherwise we consider each component
separately to get co-isometric extensions $\tilde X_j \in \Be(\Ha \oplus
\Ka_i)$ for $j$ lying in the $i$-th component of $G$. Then we denote $\Ka :=
\Ha \oplus \bigoplus_{i=1}^k \Ka_i$, where $k$ is the number of components of
$G$. We define $X_j := \tilde X_j \oplus \bigoplus_{l\neq i} I_{\Ka_l}$, where
again, $j$ lies in the $i$-th component of $G$. One easily checks that the
$X_j$'s are co-isometries that commute according to $G$.

We will proceed by induction on $n$. If $n = 1$, both statements are true by
the results about a single contraction. For the induction step assume without
loss of generality that $n$ is a vertex with only one neighbor $n-1$. Let
$\tilde G$ be the graph obtained from $G$ by deleting the vertex $n$ and
the edge $(n-1,n)$. Then $\tilde G$ is a tree and so we can apply the induction
hypothesis to $\tilde G$ and $A_1, \dots, A_{n-1}$ to get co-isometric
extensions $\tilde X_1, \dots, \tilde X_{n-1} \in \Be(\tilde \Ka)$ of $A_1,
\dots, A_{n-1}$ that commute according to $\tilde G$. Now we apply Lemma \ref
{wcl} to $A_{n-1}, A_n$ and $\tilde X_{n-1}$ to get commuting co-isometries
$X_{n-1}, X_n \in \Be(\Ka)$ that extend $\tilde X_{n-1}, A_n$, respectively.
Since $\tilde X_{n-1}$ is a co-isometry, we have $X_{n-1} = \tilde X_{n-1}
\oplus X$ and we let $X_j := \tilde X_j \oplus I_{\Ka \ominus \tilde \Ka} \in
\Be(\Ka)$, for $j = 1, 2, \dots, n-2$. It is now easy to verify that $X_1,
\dots, X_n$ commute according to $G$.

For the statement involving isometries one can verify that the induction
argument works and produces unitaries. \qed


\begin{theo} \label{main}
Let $G$ be an acyclic graph on $n$ vertices $\{1, 2, \dots, n\}$. Then for any
$n$-tuple of contractions $A_1, A_2,\dots, A_n$ on a Hilbert space $\Ha$ that
commute according to $G$, there exists an $n$-tuple of unitaries $U_1, U_2,
\dots, U_n$ on a Hilbert space $\Ka$ that commute according to $G$ and such
that
$$ P U_{j_1} U_{j_2} \dots U_{j_k} \zuz \Ha = A_{j_1} A_{j_2} \dots A_{j_k}, 
\eqno (D)$$
for all $k \in \en$, $j_i \in \{1, 2, \dots, n\}$, $1 \le i \le k$. Here
$P: \Ka \to \Ha$ is the orthogonal projection. \\ 
Conversely, if $G$ contains a cycle, there exists an $n$-tuple of contractions
that commute according to $G$ with no $n$-tuple of unitaries dilating them that
also commute according to $G$.
\end{theo} 

\pf To prove the first statement apply Lemma \ref {ind} to the contractions
$A_1^\ast, \dots, A_n^\ast$ --- they commute according to $G$. We obtain an
$n$-tuple of co-isometries $X_1, \dots, X_n \in \Be(\Ka_1)$ commuting according
to $G$ and such that $X_j$ extends $A_j^\ast$ for all $j$. In other words, we
have
$$ X_j = \begin{pmatrix} A_j^\ast & \ast \\
                                0 & \ast
         \end{pmatrix}, \hskip .5cm \hbox { so that } \hskip .5cm
  X_j^\ast = \begin{pmatrix}  A_j & 0 \\
                             \ast & \ast
             \end{pmatrix}.
$$
Now we apply the second statement of Lemma \ref {ind} to the $n$-tuple of
isometries $X_1^\ast, \dots, X_n^\ast$. We obtain unitaries $U_1, \dots, U_n
\in \Be(\Ka)$ that commute according to $G$ and extend the $X_j^\ast$'s, that
is,
$$ U_j = \begin{pmatrix} A_j &    0 & \ast \\
                        \ast & \ast & \ast \\
                           0 &    0 & \ast
         \end{pmatrix}, \hskip .5cm \hbox {for } j = 1, 2, \dots, n.
$$
An easy computation with matrices reveals that the $U_j$'s satisfy the
condition $(D)$. Thus we are done with the first part.

For the converse, we may assume that $G$ is a cycle. Indeed, otherwise we
extend the example below by taking the remaining $A_i$'s equal to the identity
operator. So, without loss of generality, the edges of $G$ are $(i,i+1)$, for
$i = 1, \dots, n-1$ and $(n,1)$ with $n \ge 3$. Let $\Ha := \Ha_0 \oplus \Ha_0$
where $\Ha_0$ is (an at least two-dimensional) Hilbert space. Let
$$ A_j := \begin{pmatrix} 0   & 0 \\
                          B_j & 0
          \end{pmatrix},
$$
where $B_j \in \Be(\Ha_0)$ is $I_{\Ha_0}$ for $j \neq 1, n$, $B_1 = R$, $B_n =
T$ and $R, T$ are non-commuting unitaries. Then $A_iA_j = 0$ for all $i,j$ and
so the $n$-tuple is commuting (and thus commuting according to $G$). It is
well-known (and easy to check) that the minimal unitary dilation of $B_2$ is
the ``inflated'' bilateral shift $\es$ on $\Ka_0 := \bigoplus_{k \in \ze}
\Ha_0$. Let us describe this in more detail. We identify the space $\Ha$ with
the sum of the zeroth and the first copy of $\Ha_0$ in $\Ka_0$, that is, $(x,y)
\in \Ha$ is equal to $\bigoplus_{k \in \ze} x_k$, with $x_0 = x$, $x_1 = y$ and
$x_k = 0$, for $k \neq 0, 1$. Then we have 
$$\es \bigoplus\nolimits_{k \in \ze} x_k = \bigoplus\nolimits_{k \in \ze} x_{k-1}.$$
Assume that there exist unitaries $U_1, \dots, U_n \in \Be(\Ka)$ dilating the
$n$-tuple and that they commute according to $G$. Then $U_2 \cong \es \oplus
\tilde U_2$ and $U_3$ has a corresponding decomposition
$$ U_3 = \begin {pmatrix} Y_{11} & Y_{12} \\
                          Y_{21} & Y_{22}
         \end {pmatrix}
       \begin {matrix} \dots & \Ka_0 \\
                       \dots & \Ka \ominus \Ka_0.
         \end {matrix}
$$
Since $U_2$ and $U_3$ commute, it follows that $\es$ commutes with $Y_{11}$
and, consequently, $Y_{11}$ is a Laurent operator on $\Ka_0 = \bigoplus_{k \in
\ze} \Ha_0$. Thus $Y_{11}$ has constant diagonals. Hence, it has $I_{\Ha_0}$ on
the diagonal just below the main diagonal. Since $U_3$ is a unitary, $Y_{11}$
is a contraction and so all the other entries of $Y_{11}$ are zero. Thus
$Y_{11} = \es$, and so it is a unitary. Since $U_3$ is a unitary, $Y_{12}$ and
$Y_{21}$ must be both zero and so $U_3$ reduces $\Ka_0$, i.e., $U_3 = \es
\oplus \tilde U_3$.

Applying the same argument repeatedly implies $U_4 = \es \oplus \tilde U_4,
 \dots, U_{n-1} = \es \oplus \tilde U_{n-1}$. For $U_1$ and $U_n$ we get
$U_1 = \re \oplus \tilde U_1$ and $U_n = \et \oplus \tilde U_n$, where $\re,
\et \in \Be(\Ka_0)$ also have constant diagonal, more precisely,
$$\re \bigoplus\nolimits_{k \in \ze} x_k = \bigoplus\nolimits_{k \in \ze}
Rx_{k-1}, \hskip .5cm \et \bigoplus\nolimits_{k \in \ze} x_k =
\ \bigoplus\nolimits_{k \in \ze} Tx_{k-1}.$$
Now we are done. Indeed, $R, T$ do not commute, hence $\re, \et$ do not commute
and so $U_1, U_n$ also do not commute, a contradiction. \qed



\begin{rem}
In the construction of the example for the converse part, we modified the
original idea of Parrott \cite {Par}. G.~Pisier pointed out to us that in his
book \cite {Pi2} there is another related example.
\end{rem}

\begin{rem}
The above theorem holds for infinite graphs.
\end{rem}

\pf Since cycles are always finite, the second part can be proven in exactly
the same manner.

To prove the first part, we use the same strategy as in the finite case. We
first observe that we may assume that $G$ is connected by the same argument.
Thus we only need to prove an infinite version of Lemma \ref {ind}. We will use
Zorn's lemma. Fix a family of contractions $\{A_v\}_{v \in V(G)} \subset
\Be(\Ha)$ that commute according to $G$. The partially ordered set $\es$ is the
set of all pairs
$$ p = ( \tilde G, \{X_v\}_{v \in V(\tilde G)} ), $$
such that $\tilde G$ is a connected subgraph of $G$, $X_v \in \Be(\Ka_p)$ is a
co-isometric extension of $A_v$ for every vertex $v$ of $\tilde G$ and
the operators $\{X_v\}_{v \in V(\tilde G)}$ commute according to $\tilde G$.
The partial order on $\es$ is given by $(\tilde G, \{X_v\}_{v \in V(\tilde G)})
\prec (G', \{Y_v\}_{v \in V(G')}) $, if and only if $\tilde G$ is a proper
subgraph of $G'$ and $Y_v$ extends $X_v$ for every $v \in V(\tilde G)$. Let 
$\{p_\lambda\}_{\lambda \in \Lambda}$ be an arbitrary chain in $\es$, where
$ p_\lambda = (G_\lambda, \{X_v^\lambda\}_{v \in V(G_\lambda)})$ and
$X_v^\lambda \in \Be(\Ka_\lambda)$. Consider the pair $p := (\tilde G,
\{X_v\}_{v \in V(\tilde G)})$, where $\tilde G = \bigcup_\lambda G_\lambda$ and
$X_v$ is the common extension of $\{X_v^\lambda\}_{\lambda \in \Lambda}$, for
each $v \in V(\tilde G)$. Since the $X_v^\lambda$'s are co-isometries, so are
the $X_v$'s. Thus the pair $(\tilde G, \{X_v\}_{v \in V(\tilde G)})$ is an
upper bound of our chain. Hence every chain in $\es$ has an upper bound and so
Zorn's lemma guarantees the existence of a maximal element $p_{\max} =
(G_{\max}, \{X_v^{\max}\}_{v \in G_{\max}})$ in $\es$. We claim that $G_{\max}
= G$ (which will complete the proof). Indeed, if $G_{\max} \subsetneq G$, then
we can find a vertex $w \in V(G)$ that is an lies on an edge $e$ whose other
end-point lies in $V(G_{\max})$. One can now find a majorant of $p_{\max}$ in
the same manner as the induction step in Lemma \ref {ind} was proved (the
corresponding graph is obtained from $G_{\max}$ by adding $w$ to the set of
vertices and $e$ to the set of edges).

The infinite version of the `isometric' part of Lemma \ref {ind} can be proved
by repeating the above argument with the obvious changes. \qed


\begin{rem}
Using the main result, we can understand the representation theory of some
universal operator algebras. Let $G$ be a graph without a cycle and let us
denote by $\A_G$ the universal operator algebra for unitaries commuting
according to $G$. That is, $\A_G$ is the norm-closed operator algebra generated
by the unitaries $\{\bigoplus_\lambda U_v^\lambda\}_{v \in V(G)}$, where $\{
\{U_v^\lambda\}_{v \in V(G)} \}_{\lambda \in \Lambda}$ are all the families of
unitaries that commute according to $G$. By definition, every contractive
representation of $\A_G$ is of the form $\bigoplus_\lambda U_v^\lambda \mapsto
A_v$, where $\{A_v\}_{v \in V(G)}$ is a family of contractions commuting
according to $G$. Conversely, if $\{A_v\}_{v \in V(G)}$ is an arbitrary family
of contractions that commute according to $G$, then we can dilate it to
unitaries $\{U_v\}_{v \in V(G)}$ commuting according to $G$, so the map $U_v
\mapsto A_v$, for all $v \in V(G)$ extends to a contractive homomorphism from
the operator algebra $\A$ generated by $\{U_v\}_{v \in V(G)}$ to the operator
algebra generated by $\{A_v\}_{v \in V(G)}$. Composing with the contractive
homomorphism from $\A_G$ to $\A$ given by $\bigoplus_\lambda U_v^\lambda
\mapsto U_v$ we get a contractive representation of $\A_G$. In fact, the
argument above shows that the representation is completely contractive. Hence
every contractive representation of $\A_G$ is completely contractive, and thus,
by Arveson's dilation theorem (Corollary 7.7 in \cite {Pau}), can be dilated to
a $\ast$-representation of $\C^\ast_G$ --- the universal C$^\ast$-algebra for
unitaries commuting according to $G$.
\end{rem}

It is possible to obtain an extension of Lemma \ref {ind} which we will now
describe. To that end, we need to explain the following natural graph
construction. Suppose we are given $n$ graphs $G_1, \dots, G_n$ (with mutually
disjoint sets of vertices), a graph $G$ on vertices $\{1,\dots, n\}$, and for
each edge $e = (k,l) \in E(G)$ a pair of vertices $v_1^e, v_2^e$ with $v_1^e
\in G_k$ and $v_2^e \in G_l$. We can then construct a graph $\hat G$ given by
$V(\hat G) = \bigcup_{i=1}^n V(G_i)$ and $E(\hat G) = \bigcup_{i=1}^n E(G_i)\ 
\cup \{(v_1^e, v_2^e);\ e \in E(G)\}$. In words, $\hat G$ is obtained by
joining the $G_i$'s together according to $G$ via edges joining the vertices
$v_1^e$'s to the corresponding $v_2^e$'s.

\begin{prop} \label {ext}
Suppose we are given graphs $\{G_i\}_{i=1}^n$ and $G$, vertices $v_1^e$,
$v_2^e$ for each $e \in E(G)$ as in the previous paragraph and contractions
$\{\{A_v^i\}_{v \in V(G_i)}\}_{i=1}^n \subset \Be(\Ha)$ that commute according
to $\hat G$. Suppose further that the graph $G$ is acyclic and that for each
fixed $i$, there exists co-isometries $\{X_v^i\}_{v \in V(G_i)} \subset
\Be(\Ka_i)$ that commute according to $G_i$ and such that $X_v^i$ extends
$A_v^i$, for each $v \in V(G_i)$. Then there exists a Hilbert space $\Ka$
containing $\Ha$, and co-isometries $\{\{\hat X_v^i\}_{v \in V(G_i)}\}_{i=1}^n
\subset \Be(\Ka)$ that commute according to $\hat G$ and such that the
$\hat X_v^i$'s extend the $A_v^i$'s.
\end{prop}
 
This result is proved by induction analogously to the proof of Lemma \ref
{ind}. The role of Lemma \ref {wcl} is played by the following fact. Given
$G_i$'s, $G$, $A_v^i$'s and $v_i^e$'s as in the above proposition with $G$
being the (unique) connected graph on two vertices and co-isometric extensions
$\{\tilde X_v^1\}_{v \in V(G_1)} \subset \Be(\tilde \Ka)$ of $\{\tilde
A_v^1\}_{v \in G_1}$ that commute according to $G_1$, there exist $\Ka \supset
\tilde \Ka$ and $X_v^i$'s as in the proposition and such that $X_v^1$ extends
$\tilde X_v^1$ for each $v \in V(G_1)$. This fact is deduced from Lemma \ref
{wcl} in essentially the same way as the induction step of the proof of Lemma
\ref {ind}.
\medskip
 
The proposition above is more general than Lemma \ref {ind} as the following
example shows. Let $G$ be the graph on vertices $\{1, 2\}$ with one edge $e$
(that joins them), $G_1$ be the cycle of length $3$ on vertices $\{w_1, w_2,
w_3\}$, and $G_2$ be the graph with one vertex $w_4$. Suppose that $v_1^e =
w_1$ and $v_2^e = w_4$ and that we are given contractions $\{A_{w_i}\}_{i=1}^4$
that commute according to $\hat G$ and such that $\{A_{w_i}\}_{i=1}^3$ possess
commuting co-isometric extensions. Then, by the proposition,
$\{A_{w_i}\}_{i=1}^4$ have co-isometric extensions which commute according to
$\hat G$. This does not follow from Lemma \ref {ind} since $\hat G$ contains a
cycle. 

\begin{rem}
One could also prove a unitary-dilation version of Proposition \ref {ext} ---
in both the statement and the proof, one would replace each occurrence of
`co-isometric extension(s)' by `unitary dilation(s)'. The downside is that one
would not be able to conclude that the property $(D)$ holds, even if it was
assumed to hold on each $G_i$ separately. This is, indeed, the reason we proved
the main theorem in two steps using co-isometric extensions.
\end{rem}
\medskip

\section* {Acknowledgements}

I would like to thank Nik Weaver for advice, John \Mcc Carthy for asking a
question that eventually led to this article and Michael Jury for pointing out
some references. In addition, I would like to thank all three of them for
helpful discussions.
\medskip

Also, I would like to thank the referee for helpful suggestions on the style
of the article and for the idea to formulate and prove Proposition \ref {ext}
(in the special case of $G$ being the connected graph on two vertices).

\end{document}